\documentclass[conference]{IEEEtran}
\IEEEoverridecommandlockouts
\usepackage{graphicx}
\usepackage{array}
\usepackage{textcomp}
\usepackage{cite}
\usepackage{url}
\usepackage{amsmath}
\usepackage{multirow}
\usepackage{bm}
\usepackage{verbatim}
\usepackage{array}
\newcolumntype{L}[1]{>{\raggedright\let\newline\\\arraybackslash\hspace{0pt}}m{#1}}
\newcolumntype{C}[1]{>{\centering\let\newline\\\arraybackslash\hspace{0pt}}m{#1}}
\newcolumntype{R}[1]{>{\raggedleft\let\newline\\\arraybackslash\hspace{0pt}}m{#1}}

\ifCLASSINFOpdf
\else
\fi
\begin{document}
\title{Optimal Allocation of Static Var Compensator via Mixed Integer Conic Programming}

\author{\IEEEauthorblockN{Xiaohu Zhang\IEEEauthorrefmark{1}$^{,}$\IEEEauthorrefmark{4},
Di Shi\IEEEauthorrefmark{1},
Zhiwei Wang\IEEEauthorrefmark{1}, 
Junhui Huang\IEEEauthorrefmark{2},
Xu Wang\IEEEauthorrefmark{2},
Guodong Liu\IEEEauthorrefmark{3} and
Kevin Tomsovic\IEEEauthorrefmark{4}}
\IEEEauthorblockA{\IEEEauthorrefmark{1}GEIRI North America}
\IEEEauthorblockA{\IEEEauthorrefmark{2}State Grid Jiangsu Electric Power Company}
\IEEEauthorblockA{\IEEEauthorrefmark{3}Oak Ridge National Laboratory}
\IEEEauthorblockA{\IEEEauthorrefmark{4}The Univerisity of Tennessee, Knoxville\\
	Email:xzhang46@vols.utk.edu}
\thanks{This work was supported in part by ARPAe (Advanced Research Projects
			Agency Energy), and in part by the Engineering Research Center Program
			of the National Science Foundation and the Department of Energy under
			NSF Award Number EEC-1041877 and the CURENT Industry Partnership Program.}}


\maketitle

\begin{abstract}
Shunt FACTS devices, such as, a Static Var Compensator (SVC), are capable of providing local reactive power compensation. They are widely used in the network to reduce the real power loss and improve the voltage profile. This paper proposes a planning model based on mixed integer conic programming (MICP) to optimally allocate SVCs in the transmission network considering load uncertainty. The load uncertainties are represented by a number of scenarios. Reformulation and linearization techniques are utilized to transform the original non-convex model into a convex second order cone programming (SOCP) model. Numerical case studies based on the IEEE 30-bus system demonstrate the effectiveness of the proposed planning model.  

\end{abstract}


% \fi
\IEEEpeerreviewmaketitle

\section*{Nomenclature}
\subsection*{Indices and Sets}
\addcontentsline{toc}{section}{Nomenclature}
\begin{IEEEdescription}[\IEEEusemathlabelsep\IEEEsetlabelwidth{$V_1,V_2,V_3$}]
	\item[$i, \ j$] Index of buses.
	\item[$k$] Index of transmission elements.
	\item[$n,m$] Index of generators and loads.
	\item[$s$] Index of scenarios.
	\item[$c$] Index of independent loop.
	\item[$\mathcal{B}$] Set of buses. 
	\item[$\mathcal{G},\mathcal{D}$] Set of generators and loads.
	\item[$\mathcal{G}_i,\mathcal{D}_i$] Set of generators and loads located at bus $i$.
	\item[$\Omega_{k}$] Set of transmission elements.
	\item[$\Omega_{s}$] Set of scenarios.
	\item[$\Omega_{c}$] Set of independent loops.
\end{IEEEdescription}

\subsection*{Variables}
\addcontentsline{toc}{section}{Nomenclature}
\begin{IEEEdescription}[\IEEEusemathlabelsep\IEEEsetlabelwidth{$V_1,V_2,V_3$}]
	\item[$P^g_{ns},Q^g_{ns}$] Active and reactive power generation of generator $n$ for scenario $s$.
	\item[$P_{ks}^r,Q_{ks}^r$] Active and reactive power flow at receiving end of branch $k$ for scenario $s$.
	\item[$P_{ks}^l,Q_{ks}^l$] Active and reactive power loss on branch $k$ for scenario $s$.
	\item[$b^v_{is}$] Susceptance of SVC at bus $i$ of scenario $s$.
	\item[$V_{is},\theta_{is}$] Voltage magnitude and angle of bus $i$ for scenario $s$.
	\item[$\delta_{k}$] Binary variable associated with installing an SVC at bus $i$.
\end{IEEEdescription} 

\subsection*{Parameters}
\addcontentsline{toc}{section}{Nomenclature}
\begin{IEEEdescription}[\IEEEusemathlabelsep\IEEEsetlabelwidth{$V_1,V_2,V_3,V_4$}]
	\item[$P_{n}^{g,\min},P_{n}^{g,\max}$] Minimum and maximum active power output of generator $n$.
	\item[$Q_{n}^{g,\min},Q_{n}^{g,\max}$] Minimum and maximum reactive power output of generator $n$.
	\item[$P_{ms}^d,Q_{ms}^d$]  Active and reactive power consumption of demand $m$ for scenario $s$.
%	\item[$S_{k}^{\max}$] Thermal limit of branch  $k$.
	\item[$V_i^{\min},V_i^{\max}$] Minimum and maximum voltage magnitude at bus $i$.
%	\item[$b^{v,\min}_{i},b^{v,\max}_{i}$] Minimum and maximum output susceptance of SVC at bus $i$.
\end{IEEEdescription} 
Other symbols are defined as required in the text.
\section{Introduction}
\label{introduction}
\IEEEPARstart{S}{tatic} Var Compensator (SVC) is one type of shunt Flexible AC Transmission Systems (FACTS) devices, which can continuously provide needed reactive power to the system \cite{mybibb:FACTS1}. The installation of SVCs at one or more appropriate locations in the network can reduce the real power loss and enhance the system transfer capability while maintaining a smooth voltage profile \cite{mybibb:FACTS2,mybibb:GA_FACT_market}. Thus, there is great interest in using SVC to improve the utilization of the existing network with the SVC allocation a concern for the system operators \cite{mybibb:GAFACTS,mybibb:svc_loadability_2013}.      

Given the nonlinear and non-convex characteristics of the load flow equations, the proper allocation of an SVC in the transmission network is a complicated task. Various heuristic approaches, such as, genetic algorithm (GA) \cite{mybibb:GAFACTS,mybibb:GA_FACT_market}, particle swarm optimization (PSO) \cite{mybibb:PSO_SQP_FACTS} and differential evolution (DE) \cite{mybibb:svc_de} have been proposed to optimally place the SVC. These techniques have the advantage of straightforward implementation but they provide no information regarding the solution quality. In \cite{mybibb:svc_sens_1}, a reactive power spot price index (QPSI) is developed to determine the best locations of an SVC. The QPSI is a weighted index at each bus under different operating conditions, including base case and critical contingencies. The authors in \cite{mybibb:svc_v_index} propose a method called the extended voltage phasors approach (EVPA) to determine the SVC locations in order to improve the voltage stability. The best SVC locations are identified based on the voltage stability index.

In \cite{mybibb:OPTSVC}, SVC allocation is formulated as a mixed integer nonlinear programming (MINLP) and solved by Benders Decomposition (BD). The objective is to maximize the system loadability considering a multi-scenario framework, including base case and several critical contingencies. A multi-start approach is embedded in the BD to avoid local optimum, and thus, the computational effort is very high. In \cite{mybibb:TCSCLFB,mybibb:svc_LFB}, the line flow based (LFB) equations are utilized to locate Thyristor Controlled Series Compensator (TCSC) and SVC. The problem is originally an MINLP model and reformulated into mixed integer linear programming (MILP) by replacing one variable in the quadratic terms with its limit. Nevertheless, the limits of some variables such as active and reactive power loss on the transmission line cannot be determined a priori. In addition, the phase angle constraints are neglected and only one load pattern is considered in the planning model. Therefore, the results obtained are suitable for preliminary analysis. 

This paper presents an optimization procedure to allocate SVC in a transmission network via mixed integer conic programming (MICP). We also leverage the LFB equations but with  improvement in several aspects. First, the quadratic loss terms are represented by second order cone (SOC) constraints and there is no need to ``guess" their limits beforehand; second, the nonlinear and non-convex term introduced by the variable susceptance of SVC is exactly reformulated into linear constraints, which maintains the SOCP format of the optimization model; third, the phase angle constraints in the meshed network is included in the model based on some standard approximations in the power flow equations; fourth, the load uncertainty is taken into account for a number of scenarios. Uncertainty considerations are increasingly important with the massive integration of the renewable energy in the power system.           

The remaining sections are organized as follows. In Section \ref{LFB}, the LFB equations in a general network are derived. The detailed formulation of the optimization model as well as the reformulation and linearization techniques are presented in Section \ref{formulation}. In Section \ref{case_study}, case studies on the IEEE 30-bus system are given. Finally, conclusions are provided in Section \ref{conclusion}. 
\section{Line Flow Based Equations in General Network}
\label{LFB}
The line flow based (LFB) equations proposed in \cite{mybibb:LFB} is originally developed for analyzing power flow in a radial distribution network. In \cite{mybibb:TEP_Joshua,mybibb:uc_lfb,mybibb:lft_kth2,mybibb:steven1}, researchers extend these equations to a transmission network and conduct various studies. For completeness, the mathematical model of LFB equations are derived in this section. 

Fig. \ref{ac_line} shows the $\pi$ model of the transmission line $k$ with bus $i$ as sending end and bus $j$ as receiving end. $r_k$, $x_k$ and $b_k$ are the resistance, reactance and shunt susceptance of line $k$. Consider the power flow direction as depicted in Fig. \ref{ac_line}, the active and reactive power balance at each bus can be expressed as follows:
\begin{figure}[!htb]
	%\captionsetup{font=footnotesize}
	\centering
	\includegraphics[width=0.3\textwidth]{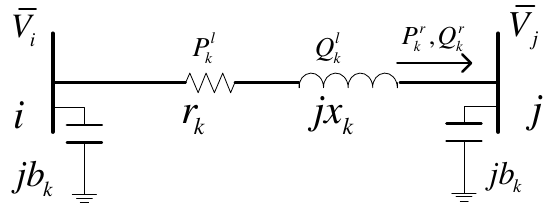}
	\caption{Static model of an AC line.}
	\label{ac_line}
\end{figure}

\begin{align}
&\sum_{n \in \mathcal{G}_i}P^g_{n}-\sum_{m \in \mathcal{D}_i}P^d_{m}=\sum_{k \in \Omega_k}\bm{M_f}(i,k)P_k^r+\sum_{k \in \Omega_k}\bm{M_l}(i,k)P^l_k  \label{real_bal}  \\
&\sum_{n \in \mathcal{G}_i}Q^g_{n}-\sum_{m \in \mathcal{D}_i}Q^d_{m}+B_iV_i^2    \nonumber   \\
&\ \ \ \ \ \ \ \ \ \ \ \ \ \ \ =\sum_{k \in \Omega_k}\bm{M_f}(i,k)Q_k^r+\sum_{k \in \Omega_k}\bm{M_l}(i,k)Q^l_k   \label{reactive_bal} 
\end{align}

In (\ref{real_bal}) and (\ref{reactive_bal}), $\bm{M_f}$ and $\bm{M_l}$ are the incidence matrix for the receiving end line flow and line losses. $\bm{M_f}(i,k)$ is 1 if bus $i$ is the sending end of line $k$, and -1 if bus $i$ is the receiving end of line $k$ and zero if neither terminal of line $k$ is connected to bus $i$. $\bm{M_l}(i,k)$ is 1 if bus $i$ is the sending end of line $k$ and zero for all the other elements. Note that $B_i$ not only includes shunt compensation but also the line susceptance connected to bus $i$. The voltage drop across branch $k$ can be written as:
\begin{equation}
\frac{V_ie^{j\theta_i}}{\tau_k e^{j\theta_k^{ps}}}=V_je^{j\theta_j}+\frac{P_k^r-jQ_k^r}{V_je^{-j\theta_j}}(r_k+jx_k) \label{v_drop}
\end{equation}  

Equation (\ref{v_drop}) is valid for a general branch element, including the phase shifting transformer and a transmission line. For a normal transmission line, the tap ratio $\tau_k$ is set to be 1 and the phase shifting $\theta_k^{ps}$ is set to be zero. Taking the magnitude and imaginary part of (\ref{v_drop}), the followings are obtained:
\begin{align}
&V_i^2/\tau_k^2-V_j^2=2r_kP^r_k+2x_kQ^r_k+r_kP_k^l+x_kQ_k^l   \label{v_drop_m}  \\
&V_iV_j\sin(\theta_i-\theta_j-\theta^{ps}_k)=\tau_k(P^r_kx_k-Q^r_kr_k)  \label{v_drop_im}
\end{align}
Assuming the bus voltage magnitude is close to 1 and the phase angle difference across a branch is small, the equation (\ref{v_drop_im}) can be approximated as:
\begin{equation}
\theta_i-\theta_j \approx \tau_k(P^r_kx_k-Q^r_kr_k)+\theta_k^{ps}    \label{approx_angle}
\end{equation}

To manage the angle constraints in the meshed network, the following equation is introduced \cite{mybibb:lft_kth2}:
\begin{equation}
\sum_{k \in \Omega_k}\bm{C_l}(c,k)(\tau_k(P^r_kx_k-Q^r_kr_k)+\theta_k^{ps}) \approx 0  \label{loop}
\end{equation}
Equation (\ref{loop}) expresses the fact that the phase angle difference across each independent loop of a graph is zero. $\bm{C_l}$ is the incidence matrix corresponding to all the independent loops in the network. With a predefined loop direction, $\bm{C_l}(c,k)$ is 1 if branch $k$ is with the same direction as loop $c$, and -1 if branch $k$ is with the opposite direction of loop $c$ and zero if branch $k$ is not in loop $c$. Due to the approximation used in (\ref{approx_angle}), equation (\ref{loop}) may not hold in the strict sense, we slightly relax the constraint by introducing a small allowable error $\epsilon_{\theta}$ ($\pi/360$ or $0.5 \deg$) and obtain:
\begin{equation}
-\epsilon_{\theta} \le \sum_{k \in \Omega_k}\bm{C_l}(c,k)(\tau_k(P^r_kx_k-Q^r_kr_k)+\theta_k^{ps}) \le \epsilon_{\theta}  \label{relax_loop}
\end{equation}

The active and reactive power loss on each branch can be described as:
\begin{align}
& P^l_k=\frac{(P^{r}_k)^2+(Q_k^r)^2}{V_j^2}r_k  \label{soc1} \\
&P^l_kx_k=Q^l_kr_k  \label{loss_relation}   
\end{align}
With an auxiliary variable $P^{l,aux}_k$, constraint (\ref{soc1}) can be transformed into a rotated quadratic cone by relaxing the equality into inequality \cite{mybibb:lft_kth2}:
\begin{align}
&2P^{l,aux}_kV_j^2 \ge (P^{r}_k)^2+(Q_k^r)^2  \label{soc2}   \\
&P^l_k=2r_kP^{l,aux}_k      \label{aux_c} 
\end{align} 

For the accuracy of the LFB equations in optimal power flow (OPF) problem, we refer readers to \cite{mybibb:TEP_Joshua,mybibb:lft_kth2}, which compare the LFB based OPF results with standard NLP based OPF results for different IEEE test systems.

\section{Problem Formulation}
\label{formulation}
In this section, the complete optimization model SVC allocation is first presented, and then reformulation and linearization techniques are leveraged to transform the non-convex model into a convex SOCP model. 
\subsection{Optimization Model}
The SVC can provide local reactive power compensation for the power system so it is suitable for minimizing the transmission loss and regulating voltage. The complete optimization model is given by (\ref{obj})-(\ref{sum_svc}):
\begin{align}
&\min_{\Xi_{\text{OM}}} A_1 \sum_{s \in \Omega_s}\rho_s \sum_{k\in \Omega_k}P^l_{ks}+A_2 \sum_{s \in \Omega_s}\rho_s \sum_{i\in \mathcal{B}}|V^2_{is}-1|    \nonumber \\   
&\ \ \ \ \ \ \ \ +\alpha \sum_{s \in \Omega_s}\rho_s \sum_{k\in \Omega_k}P^{l,aux}_{ks}  \label{obj} \\
& \text{subject to:}    \nonumber  \\
&(\ref{real_bal}),(\ref{v_drop_m}),(\ref{relax_loop}),(\ref{loss_relation})-(\ref{aux_c}) \ \ \text{and}  \nonumber  \\
&\sum_{n \in \mathcal{G}_i}Q^g_{ns}-\sum_{m \in \mathcal{D}_i}Q^d_{ms}+B_iV_{is}^2+\delta_ib^v_{is}V_{is}^2    \nonumber   \\
&\ \ \ \ \ \ \ \ \ \ \  =\sum_{k \in \Omega_k}\bm{M_f}(i,k)Q_{ks}^r+\sum_{k \in \Omega_k}\bm{M_l}(i,k)Q^l_{ks}   \label{reactive_bal_svc} \\
&(S^{\max}_{k})^2 \ge (P^r_{ks})^2+(Q^r_{ks}+V_{js}^2b_k)^2    \label{to_thermal}  \\
&(S^{\max}_{k})^2 \ge (P^r_{ks}+P^l_{ks})^2+(Q^r_{ks}+Q^l_{ks}-V_{is}^2b_k)^2   \label{from_thermal}  \\
&(V_i^{\min})^2 \le V_{is}^2 \le (V_i^{\max})^2    \label{volt_lim}  \\
&P^{g,\min}_n \le P^{g}_{ns} \le P^{g,\max}_n    \label{gp_lim}   \\
&Q^{g,\min}_n \le Q^{g}_{ns} \le Q^{g,\max}_n    \label{gq_lim}   \\
&b^{v,\min}_i \le b_{is}^{v} \le b_i^{v,\max}    \label{svc_lim}  \\
&\sum_{i\in \mathcal{B}}\delta_i \le N_v     \label{sum_svc}
\end{align}
Constraints (\ref{obj})-(\ref{sum_svc}) hold $\forall c \in \Omega_c, s \in \Omega_s, n\in \mathcal{G}, i,j\in \mathcal{B}, k\in \Omega_k$.

The optimization variables of the SVC allocation model are those in the set $\Xi_{\text{OM}}=\{P^g_{ns},Q^g_{ns},P^r_{ks},Q^r_{ks},P^l_{ks},Q^l_{ks}, V_{is}^2,\delta_i,\\ b_{is}^v,P^{l,aux}_{ks}\}$. The objective function contains three terms that are all weighted by the probability $\rho_s$ for each scenario. The first term is to minimize the network loss, which is the summation of the loss on all the branches. The second term is to minimize the voltage deviation to improve the voltage profile. $A_1$ and $A_2$ are the weighting factor for these two objectives. The third term is the penalty term for driving the rotated cone constraint (\ref{soc2}) to be binding at the optimal solution \cite{mybibb:lfb_kth3}. $\alpha$ is a scaling factor for the penalty term which should be carefully selected to ensure that the penalty term is only a small fraction of the objective function and does not adversely affect the allocation strategy.   

Constraint (\ref{reactive_bal_svc}) represents the reactive power balance at each bus. A binary variable $\delta_i$ is introduced to flag the installation of an SVC. The receiving and sending end thermal limits of a branch are enforced by constraints (\ref{to_thermal})-(\ref{from_thermal}). $S^{\max}_k$ is the thermal limit for branch $k$. Note that these two constraints are also the second order cone constraints. The physical limits for voltage magnitude, active and reactive power generation and SVC output susceptance are indicated by constraints (\ref{volt_lim})-(\ref{svc_lim}). Constraint (\ref{sum_svc}) limits the number of SVCs that can be installed in the system to $N_v$. Note that we assume a predetermined number of SVC are ready for installation so the capital cost of SVC is a sunk cost and not included in the objective. Nevertheless, the installation cost of SVC can be easily embedded in the objective function with an appropriate scaling factor.
         
\subsection{Reformulation and Linearization}  
As observed from the optimization model, there exists two terms which are not convex and need to be reformulated and linearized. In the objective function, the absolute value can be cast into an LP by introducing two positive slack variables and one additional constraint:
\begin{align}
&V_i^2-1+s_{i,1}-s_{i,2}=0   \label{abs_cons}  \\
&s_{i,1} \ge 0, s_{i,2} \ge 0   \label{slack_con}
\end{align}  
Then the second term of the objective function becomes:
\begin{equation}
A_2 \sum_{s \in \Omega_s}\rho_s \sum_{i\in \mathcal{B}}|V^2_{is}-1| \rightarrow A_2 \sum_{s \in \Omega_s}\rho_s \sum_{i\in \mathcal{B}}(s_{is,1}+s_{is,2})
\end{equation}

In constraint (\ref{reactive_bal_svc}), there is a trilinear term $\delta_ib^v_iV_i^2$ which involves the product of a binary variable and two continuous variables. An approach similar to the reformulation technique proposed in \cite{mybibb:investment_naps} is used to linearize the trilinear term. We first introduce a new variable $Q^v_i$:
\begin{equation}
Q^v_i=\delta_ib^v_iV_i^2   \label{new_variable}
\end{equation}
We then multiply each side of constraint (\ref{svc_lim}) with $\delta_i$ and combine with (\ref{new_variable}) to yield:
\begin{equation}
\delta_i b^{v,\min}_i \le Q^v_i/V_i^2 \le \delta_i b^{v,\max}_i    \label{reformulation_1}
\end{equation}
Since $V_i^2$ is always larger than zero, constraint (\ref{reformulation_1}) can be written as:
\begin{equation}
\delta_i V_i^2b^{v,\min}_i \le Q^v_i \le \delta_i V_i^2b^{v,\max}_i    \label{reformulation_2}
\end{equation}
Constraint (\ref{reformulation_2}) still includes a nonlinear term which is the product between a binary variable and a continuous variable. This term can be linearized by introducing another variable $z_i=\delta_i V_i^2$:
\begin{align}
& z_ib^{v,\min}_i \le Q^v_i \le z_ib^{v,\max}_i   \label{refor_1}  \\
&\delta_i(V_{i}^{\min})^2 \le z_i \le \delta_i(V_{i}^{\max})^2 \label{refor_2}  \\
&V_i^2-(1-\delta_i)(V_{i}^{\max})^2 \le z_i \le V_i^2-(1-\delta_i)(V_{i}^{\min})^2 \label{refor_3}
\end{align}

Therefore, the trilinear term $\delta_ib_k^vV_i^2$ is linearized by three constraints from (\ref{refor_1})-(\ref{refor_3}). Finally, all the constraints in the optimization model are either linear or in conic formats, which can be solved efficiently by commercial solvers.
\section{Case Studies}
\label{case_study}
The IEEE-30 bus system is selected to test the performance of the proposed optimization procedure. It has six generators, four transformers and 37 transmission lines. The base active and reactive loads are 260 MW and 116 MVar. Additional system data can be found in MATPOWER software package \cite{mybibb:MATPOW}. The problem is modeled in YALMIP \cite{mybibb:YALMIP} and solved by MOSEK \cite{mybibb:mosek}. The computer used to perform the computational tasks has an Inter Core(TM) i5-2400M CPU @ 2.30 GHz and 4.00 GB of RAM. 

The allowable compensation range for the SVC is from 0 to 0.3 p.u. \cite{mybibb:PSO_SQP_FACTS}. We assume that every bus, except generator buses, is a candidate location to install SVC so the number of binary variables is 24. The load uncertainties are represented by 15 scenarios as given in Table \ref{data_scenario} \cite{mybibb:tcsc_ac}. Note that the actual load for each scenario is the base load multiplied by the load factor $\lambda_s$. 
\begin{table}[!htb]
	\centering
	\caption{Probability and Load Factor for Each Scenario}
	\label{data_scenario}
	\begin{tabular}{|c|c|c|c|c|c|c|c|c|}
		\hline
		Scenario&1&2&3&4&5&6&7&8   \\
		\hline
		$\rho_s$&0.02&0.14&0.04&0.02&0.14&0.04&0.02&0.14         \\
		\hline
		$\lambda_s$&1.00&0.80&0.60&1.10&0.88&0.66&1.21&0.97 \\
		\hline
		Scenario&9&10&11&12&13&14&15&   \\
		\hline
		$\rho_s$&0.04&0.02&0.14&0.04&0.02&0.14&0.04&         \\
		\hline
		$\lambda_s$&0.73&1.33&1.06&0.77&1.46&1.17&0.87& \\
		\hline
	\end{tabular}
\end{table}     

Obviously, the optimization results depend on the weighting coefficients $A_1$ and $A_2$. We therefore test four different weighting schemes as given in Table \ref{diff_case}. Based on the trial-and-error analysis, $\alpha=0.001$ is sufficient for the rotated cone constraints to be binding at the optimality for all cases.  
\begin{table}[!htb]
\centering
\caption{Difference Weighting Schemes}
\label{diff_case}
\begin{tabular}{|c|c|c|c|c|}
 \hline
 & Case 1&Case 2&Case 3&Case 4    \\
 \hline
 $A_1$&1&1&10&1     \\
 \hline
 $A_2$&0&1&1&10    \\
 \hline
\end{tabular}
\end{table} 

Table \ref{res_diff_case} provides SVC placement solutions under different values of $N_v$. In the first column, $P^l_0$ and $\Delta V_0$ are the weighted loss and voltage deviation without any SVC. The third and fourth column give the weighted loss and voltage deviation with the SVCs. The selected locations are provided in column five. Column six shows the maximum cone mismatch for constraint (\ref{soc2}) considering all the scenarios. The last column presents the computational time for each case.   

\begin{table}[!htb]
	\centering
	\caption{Allocation Strategy of SVC for Different Cases}
	\label{res_diff_case}
	\begin{tabular}{|c|c|c|c|c|c|c|}
		\hline
		&\multirow{3}{0.2cm}{$N_v$}&\multirow{3}{0.6 cm}{$\ \ P^l$ (MW)}&\multirow{3}{0.5cm}{$\Delta V$ \newline (p.u.)}&&Max&\multirow{3}{*}{\parbox{0.7cm}{\centering Time (s)}}    \\
		Case &&&&Location&cone&    \\
		&&&&&error&     \\
		\hline
		Case 1&1&2.55&1.24&21&2e-5&16.18    \\
		
		$P^l_{0}$&2&2.48&1.28&4,24&5e-5&41.71   \\
		
		(2.70 MW)&3&2.42&1.39&4,21,24&8e-6&74.16   \\
		
		$\Delta V_0$&4&2.40&1.43&4,19,21,24&2e-5&214.94    \\
		
		(1.08 p.u.)&5&2.39&1.45&4,19,21,24,26&2e-5&265.08    \\
		\hline
		Case 2&1&2.75&0.22&24&6e-6&15.07    \\
		
		$P^l_{0}$&2&2.74&0.16&19,24&7e-6&55.16   \\
		
		(2.81 MW)&3&2.66&0.15&4,19,24&2e-5&115.79   \\
		
		$\Delta V_0$&4&2.63&0.14&4,19,21,24&2e-5&308.20    \\
		
		(0.33 p.u.)&5&2.62&0.11&4,19,21,24,30&2e-5&402.86    \\
		\hline
		 Case 3&1&2.63&0.45&24&5e-5&18.50    \\
		 
		 $P^l_{0}$&2&2.56&0.44&4,24&7e-5&34.29   \\
		 
		 (2.75 MW)&3&2.53&0.45&4,21,24&9e-5&91.26   \\
		 
		 $\Delta V_0$&4&2.51&0.42&4,19,21,24&1e-5&210.63    \\
		 
		 (0.43 p.u.)&5&2.50&0.39&4,7,19,21,24&7e-6&315.51    \\
		 \hline
		 Case 4&1&2.92&0.18&24&4e-5&14.38    \\
		 
		 $P^l_{0}$&2&2.84&0.14&19,24&2e-5&55.48   \\
		 
		 (2.86 MW)&3&2.84&0.11&6,19,24&5e-5&131.08   \\
		 
		 $\Delta V_0$&4&2.81&0.10&6,19,24,30&3e-5&295.77    \\
		 
		 (0.33 p.u.)&5&2.80&0.08&4,7,19,24,30&1e-5&408.95    \\
		 \hline
	\end{tabular}
\end{table} 

As observed from Table \ref{res_diff_case}, the allocation strategy varies in different cases. Case 1 refers to a single objective optimization of the real power loss. Without SVC, the network loss is 2.70 MW. With 5 SVCs placed in the system, the loss decreases to 2.39 MW. Note that since the voltage deviation is not considered in the objective of case 1, it increases as the power loss decreases. In case 2, the weighting coefficient for the two objectives are equal. Both the power loss and the voltage deviation decrease as the number of SVCs increases. In case 3, the weighting coefficient for the power loss is 10 times of that for the voltage deviation. Compared to the objective without SVC, it can be seen that the more SVCs the network has, the more reductions of power loss can be observed. However, this is not the situation for the voltage deviation. There is a slight increase on the voltage deviation when no more than 3 SVCs are allowed to be installed in the system. When there are more than 4 SVCs placed in the network, the voltage deviation also decreases. Similar results can also be observed in case 4 when the weighting factor for the voltage deviation is higher. The power loss will increase when there is only 1 SVC in the system and it will be reduced when more than 2 SVCs are allowed to be placed in the network. For all the cases, the computational time increases as the maximum number of SVCs ($N_v$) increases.

\begin{figure}[!htb]
	%\captionsetup{font=footnotesize}
	\centering
	\includegraphics[width=0.46\textwidth]{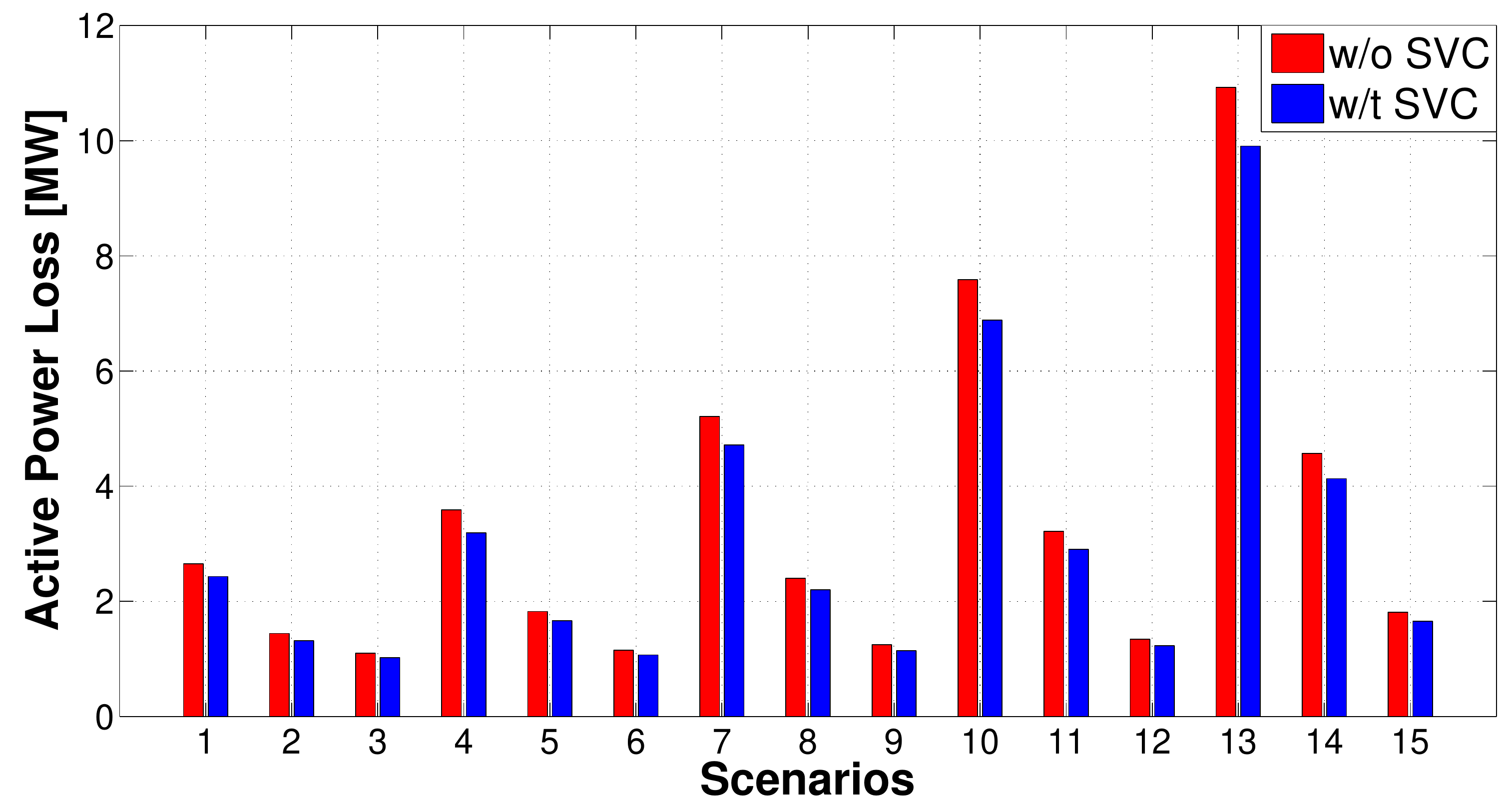}
	\caption{Network loss for each scenario in case 3.}
	\label{loss_comp}
\end{figure}

Fig. \ref{loss_comp} depicts the real power loss for each scenario in case 3. The blue bar refers to the network loss with 5 SVCs. It can be seen that the loss reduction can be observed in all scenarios and the most reduction occurs in scenario 13 in which the load demand is the highest.  

\begin{figure}[!htb]
	%\captionsetup{font=footnotesize}
	\centering
	\includegraphics[width=0.46\textwidth]{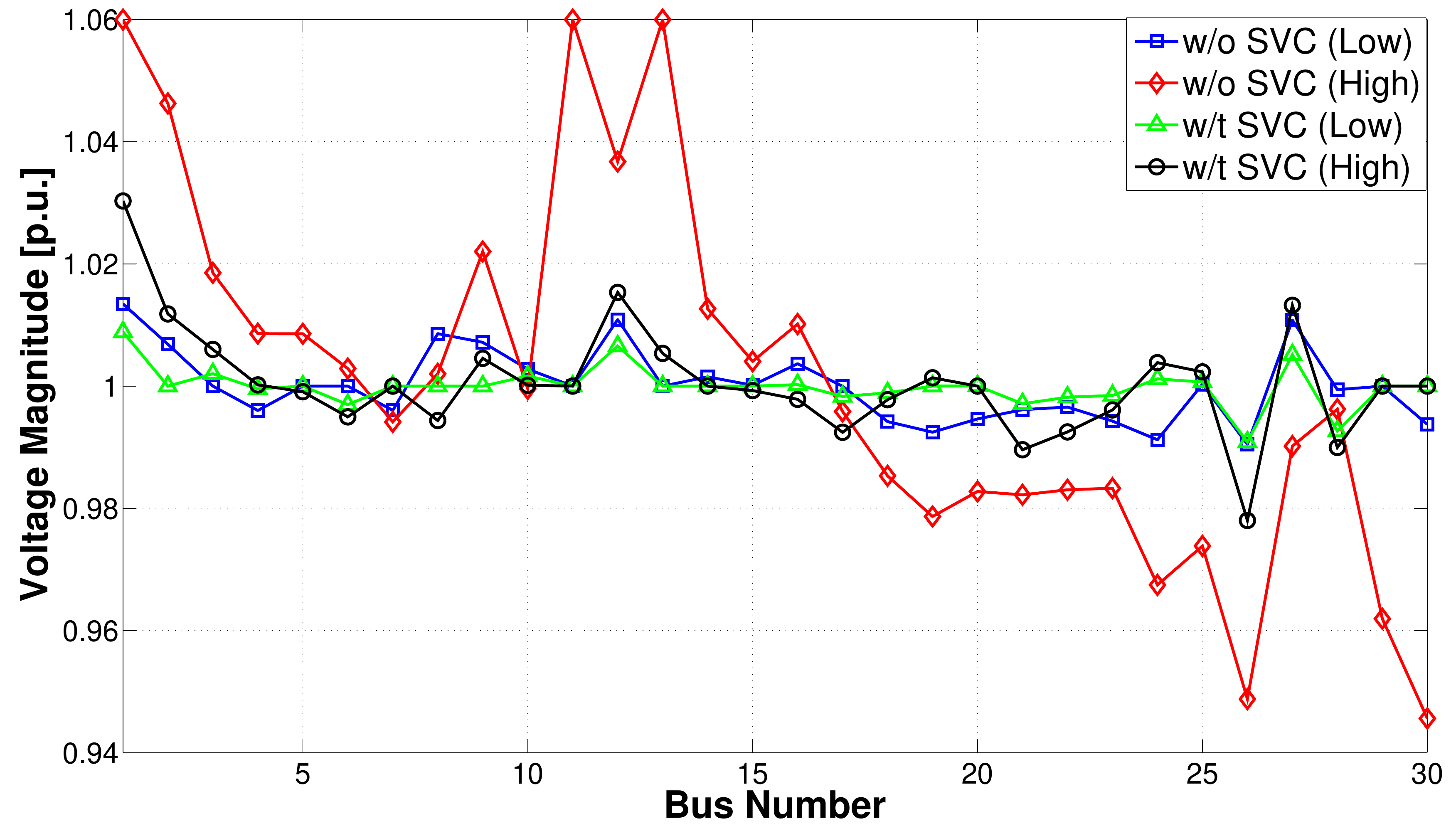}
	\caption{Voltage magnitude at each bus in case 4 for two different scenarios.}
	\label{v_comp}
\end{figure}

Fig. \ref{v_comp} shows the voltage magnitude at each bus in case 4 for the lowest load level (scenario 3) and the highest load level (scenario 13). The green and black lines refer to the voltage magnitude with 5 SVCs. In both of the two scenarios, the installation of SVCs will push the voltage magnitude close to 1 p.u. and thus, a better voltage profile is achieved.     
       
\section{Conclusions}
\label{conclusion}
This paper presents a planning model to determine the suitable locations of SVC in a transmission network to minimize the system real power loss and improve the voltage profile. The load uncertainties are modeled by a number of scenarios. The planning model is originally a large scale MINLP model which is difficult to solve. Reformulation and linearization techniques are used to convexify the planning model and transform it into MICP, which can be efficiently solved by commercial solvers. Numerical case studies based on the IEEE 30-bus system demonstrate the effectiveness of the proposed planning model. Simulations show that if several SVCs are appropriately placed in the system, the real power loss can be reduced and a smoother voltage profile can be achieved.

\bibliographystyle{IEEEtran}
\bibliography{IEEEabrv,mybibb}
%\end{thebibliography}
%\bibliographystyle{IEEEtran}
%\bibliography{IEEEabrv,mybibb}

\end{document}